\documentclass[11pt]{article}

\usepackage{colordvi,epsfig,amsmath,amssymb}

\newtheorem{teo}{Theorem}[section]
\newtheorem{lem}{Lemma}[section]

\newcommand{\ignore}[1]{}

\newcommand{\qed}{\hfill $\square$}

\begin{document}

\title{Convergence Analysis of an Inexact Feasible Interior Point Method 
       for Convex Quadratic Programming\thanks{Supported by EPSRC grant EP/I017127/1.}}

\author{Jacek Gondzio\thanks{Email: {\tt J.Gondzio@ed.ac.uk},
          URL: {\tt http://maths.ed.ac.uk/\~{}gondzio/ }}
          \vspace{1.0cm}
        \\
School of Mathematics \\
and Maxwell Institute for Mathematical Sciences \\
The University of Edinburgh \\
Mayfield Road, Edinburgh EH9 3JZ \\
United Kingdom. \\ \\ \\ 
Technical Report ERGO-2012-008\thanks{For other papers in this series see
         {\tt http://www.maths.ed.ac.uk/ERGO/}}
\vspace{+.7cm} \\
}
\date{July 25, 2012}

\begin{titlepage}
\maketitle

\begin{abstract}

In this paper we will discuss two variants of an inexact feasible 
interior point algorithm for convex quadratic programming. 
We will consider two different neighbourhoods: 
a (small) one induced by the use of the Euclidean norm which yields 
a short-step algorithm and a symmetric one induced by the use 
of the infinity norm which yields a (practical) long-step algorithm. 
Both algorithms allow for the Newton equation system to be solved 
inexactly. For both algorithms we will provide conditions for 
the level of error acceptable in the Newton equation and establish 
the worst-case complexity results.

\end{abstract}

\hspace{0.5cm}

\noindent {\it Keywords:}  
Inexact Newton Method, Interior Point Algorithms, Linear Programming, 
Quadratic Programming, Worst-case Complexity Analysis, Matrix-Free Methods.

\end{titlepage}

\pagestyle{myheadings}
\markboth{
\underline{Inexact Feasible IPMs}}
{\underline{Inexact Feasible IPMs}}
\parindent 0pt \parskip 10pt

\newpage

\setcounter{page}{2}
\section{Introduction}
\label{Intro}

It is broadly accepted that interior point methods (IPMs) provide 
very efficient solution techniques for linear and convex quadratic 
programming problems \cite{iter:G-ipmXXV,iter:Wright}.
Interior point algorithms for these classes of problems enjoy 
excellent worst-case complexity bounds: indeed the best known 
algorithms find the $\varepsilon$-accurate optimal solutions 
to the problem with $n$ variables in ${\cal O}(\sqrt{n} \log(1/\varepsilon))$ 
or ${\cal O}(n \log(1/\varepsilon))$ iterations, depending 
on how aggressive steps to optimality are allowed. 
Computational experience provides evidence that the algorithm 
which uses a more aggressive strategy (the so-called long-step 
method) solves linear and quadratic programming problems 
in a number of iterations which may be expressed 
as ${\cal O}(\log{n} \log(1/\varepsilon))$ \cite{iter:G-ipmXXV}.

The small number of iterations does not always guarantee the efficiency 
of the method because occasionally IPMs struggle with a high 
per-iteration cost of the linear algebra operations. In the most 
adverse case the cost of solving a dense optimization problem 
employing a direct linear algebra method to solve the Newton 
equation system may reach ${\cal O}(n^3)$ flops per iteration. 
The effort of a single IPM iteration is usually significantly 
lower than this upper bound. However if problems are very large 
then, although they may display reasonable sparsity features, 
the use of direct sparsity-exploiting linear algebra techniques 
may still run into trouble due to excessive memory requirements 
or unacceptably long CPU time. 
Iterative methods for linear equations such as conjugate 
gradients or other approaches from the Krylov subspace family 
may offer a viable alternative to direct methods in such cases. 

The interest in the use of iterative methods to solve the Newton 
equation system in IPMs has been growing over the last decade: 
see for example the recent survey of D'Apuzzo et al. \cite{iter:dAdSdS}.
Iterative methods offer several advantages. In particular they are 
often memory efficient and hence allow much larger problems to be solved. 
However to take full advantage of iterative methods one usually 
has to relax the accuracy requirements in the solution of the 
Newton equation system. Consequently, instead of using an {\it exact} 
Newton direction, the resulting IPM employs an {\it inexact} one. 
This opens up all sorts of theoretical and practical questions. 
We will state and answer some of these in this paper. 
In particular, we will discuss the key issue concerning an acceptable 
level of error in the inexact Newton method used by an IPM. 

The use of an inexact Newton method \cite{iter:DES} is well established 
in the context of solving nonlinear equations \cite{iter:Kelley} 
and in nonlinear optimization \cite{iter:NocedalWright}. 
Bellavia \cite{iter:Bellavia} applied an inexact interior point 
method to solve monotone nonlinear complementarity problems 
and proved global convergence and local superlinear convergence 
of the method. Several successful attempts were also made to shed 
light on the application of an inexact Newton method in IPMs for linear 
and convex quadratic programming problems. 

Freund, Jarre and Mizuno \cite{iter:FJM-inexact} and Mizuno and 
Jarre \cite{iter:MJ-inexact} extended a very popular globally 
convergent infeasible path-following method for linear programming 
of Kojima, Megiddo and Mizuno \cite{iter:KMM-Inf-PD} to accommodate 
the inexact solution of Newton systems. In particular Mizuno 
and Jarre \cite{iter:MJ-inexact} proved that an inexact variant 
of this algorithm has ${\cal O}(n^2 \log(1/\varepsilon))$ complexity. 
Baryamureeba and Steihaug \cite{iter:VB-TS} provided another 
extension of the method of Kojima et al., allowing for inexactness 
in both primal and dual Newton steps. All these analyses considered 
a general case in which no assumption was made about {\it how} 
the Newton system is solved. The only assumptions made were concerned 
with the absolute or relative error in the inexact Newton direction. 

Several authors tried to specialize their analyses using a better 
understanding of the specific iterative methods of linear algebra 
employed to solve the Newton equation system inexactly. Indeed, 
if a particular preconditioner is used in a given Krylov subspace 
method applied in this context, it is often possible to control 
the residuals in specific linear subspaces and design special 
variants of an inexact interior point algorithm. 
Al-Jeiroudi and Gondzio \cite{iter:alJG} considered the case 
in which an indefinite preconditioner based on a guess 
of basic-nonbasic partition inspired by the work of Oliveira 
and Sorensen \cite{iter:OS-pcg1} is used to precondition 
the indefinite augmented system (the reduced KKT system). 
They designed an inexact infeasible primal-dual IPM for linear 
programming and established its ${\cal O}(n^2 \log(1/\varepsilon))$ 
complexity. 
Lu, Monteiro and O'Neal \cite{iter:LMoN} analysed the case 
of quadratic programming in which the matrix of the quadratic 
objective term has a known factorization and proposed an interesting 
specialized preconditioner for the Newton system in this case. 
They showed that the resulting inexact IPM converges 
in ${\cal O}(n^2 \log(1/\varepsilon))$ iterations. 
Cafieri et al. \cite{iter:CdAdSdST} performed an analysis 
of an inexact potential reduction algorithm for convex 
quadratic programming.

In this paper we consider a {\it feasible} primal-dual path-following 
method for convex quadratic programming and analyse it in a situation 
when the solutions of the Newton systems admit a certain level of inaccuracy. 
We prove worst-case iteration complexity results for two variants 
of such a method. The first variant requires the iterates to stay 
in a small neighbourhood of the central path induced by the use 
of the Euclidean norm to control the error in the perturbed 
complementarity conditions. Such a method has the best known 
iteration complexity result ${\cal O}(\sqrt{n} \log(1/\varepsilon))$ 
and we prove that the inexact variant preserves the same complexity. 
However, this method is only of theoretical interest because its 
implementation demonstrates the behaviour predicted by the worst-case 
analysis. Therefore, this algorithm is not used in practice. 
The second variant allows the iterates to stay in a wide symmetric 
neighbourhood of the central path induced by the use of the infinity 
norm to control the error in the perturbed complementarity conditions. 
We show that this method reaches an $\varepsilon$-optimal solution 
in ${\cal O}(n \log(1/\varepsilon))$ iterations. The second 
algorithm has a practical meaning and although its worst-case 
iteration complexity is $\sqrt{n}$ times higher than the former 
one, in practice this algorithm behaves very well and provides 
the basis of an implementable method. 

To simplify the analysis and to allow the reader to concentrate 
on the essential consequences of inexactness in the Newton 
direction we will analyse {\it feasible} algorithms. 
It is worth mentioning that this does not limit the applicability 
of the analysis. Indeed, the homogeneous and self-dual emdedding 
\cite{iter:YTM-self-dual} provided an elegant tool to transform 
the linear programs into new models for which the primal and dual 
initial feasible solutions are known and therefore primal and 
dual feasibility can be maintained throughout the computations. 
The homogeneous and self-dual emdedding was initially used 
in the context of LPs and implemented by Andersen and Andersen 
\cite{iter:AA-MOSEK} in their Mosek software. Later this elegant model 
was extended to conic optimization by Andersen et al. \cite{iter:ART-conic}.

Finally, let us comment that the key motivation for this work 
is the need to better understand how much inaccuracy is admissible 
in the Newton systems and provide the foundations for cases in which 
interior point algorithms rely on iterative methods to solve 
the underlying linear algebra problems. There has recently been 
a shift of interest in the IPM community towards the application 
of iterative methods to solve the reduced KKT systems \cite{iter:dAdSdS}. 
There exists a rich body of literature (cf. Benzi et al. \cite{iter:BenziGL}) 
which deals with very similar saddle point problem arising 
in the discretisations of partial differential equations and 
we expect that the coming years will bring many interesting 
developments in iterative methods applicable to IPMs. 
In particular there is a clear need to develop alternative 
preconditioners which are compatible with the spirit of the matrix-free 
interior point method \cite{iter:G-ipmXXV,iter:G-MtxFree}.  

%

The paper is organised as follows. In Section~\ref{IPMsBack} we will 
introduce the quadratic optimization problem, define the notation 
used in the paper and point out an essential difference between 
the {\it exact} \, and {\it inexact} \, interior point methods.
In Section~\ref{Complexity} we will perform the worst-case analysis 
of two variants of an inexact feasible interior point algorithm 
for convex quadratic programming. First we will analyse the algorithm 
operating in a small neighbourhood of the central path induced 
by the 2-norm. Such a method yields the best complexity result known 
to date but it has only a theoretical importance. Next we will 
analyse the feasible algorithm operating in a symmetric neighbourhood 
of the central path induced by the infinity norm. This method has 
a practical meaning as an implementable algorithm. It provides 
a theoretical basis for the recently developed {\it matrix-free} 
variant of the interior point method \cite{iter:G-MtxFree}. 
Our analysis will follow that of Wright \cite{iter:Wright} and will 
generalize it from linear programming to quadratic programming 
and from exact to inexact method. 
In Section~\ref{PracIPM} we will briefly comment on some practical 
aspects related to the implementation of inexact interior point 
method and finally in Section~\ref{Conclusions} we will give our 
conclusions.

\section{Interior point methods: background}
\label{IPMsBack}
We are concerned in this paper with the theory of interior point 
methods for solving convex quadratic programming (QP) problems. 
We consider the following general primal-dual pair of QPs 
\begin{eqnarray}
   \begin{array}{rlcrl}
   \mbox{ {\small Primal} } & & \quad \quad \quad \quad \quad &
   \mbox{ {\small Dual} }   &     \vspace{0.3cm} \\
   \mbox{ min }   & c^T x + \frac{1}{2} x^T Q x & &
   \mbox{ max }   & b^T y - \frac{1}{2} x^T Q x \\
   \mbox{ s.t. }  &   A x = b,                  & &
   \mbox{ s.t. }  & A^T y + s - Q x = c,        \\
                  &     x \geq 0, &  &
                  & y \mbox{ free,} \ \ s \geq 0,  
   \end{array}
   \label{QPpair}
\end{eqnarray}
where $A \in {\cal R}^{m \times n}$ has full row rank $m \leq n$, 
$\, Q \in {\cal R}^{n \times n}$ 
is a positive semidefinite matrix,
$x, s, c \in {\cal R}^{n}$ and $y, b \in {\cal R}^{m}$.
Setting $Q=0$ yields the special case of the linear programming (LP) 
primal-dual pair. 

To derive a primal-dual interior point method \cite{iter:Wright} we first 
introduce the logarithmic barrier function $\mu \sum_{j=1}^{n} \log x_j$ 
to ``replace'' the inequality constraint $x \geq 0$ in the primal and 
then, by using Lagrangian duality theory, write down the first-order 
optimality conditions
\begin{eqnarray}
  \begin{array}{rcl}
     A x             & = & b, \\
     A^T y + s - Q x & = & c, \\
     X S e           & = & \mu e, \\
     (x,s)           & \geq & 0, 
  \end{array}
   \label{FOC_barrQP}
\end{eqnarray}
where $X$ and $S$ are diagonal matrices in ${\cal R}^{n \times n}$
with elements of vectors $x$ and $s$ spread across the diagonal,
respectively and $e \in {\cal R}^n$ is the vector of ones.
This system of equations has a unique solution
$(x(\mu), y(\mu), s(\mu)), \, x(\mu)>0, \, s(\mu)>0$ for any $\mu > 0$. 
The corresponding point is called a $\mu$-centre. A family of these 
points for all positive values of $\mu$ determines a continuous 
curve $\{ (x(\mu), y(\mu), s(\mu)): \mu > 0 \}$ which is called 
the primal-dual {\it central trajectory} or {\it central path}.

The first two equations in (\ref{FOC_barrQP}) are the primal and dual 
feasibility conditions, respectively. The third equation is a perturbed 
complementarity condition; the parameter $\mu$ associated with 
the logarithmic barrier function controls the level of perturbation. 
The convergence to optimality is forced by taking the barrier 
parameter $\mu$ to zero \cite{iter:Wright,iter:G-ipmXXV}.  

In this paper we will assume that we work with the {\it feasible} 
interior point method and therefore all iterates $(x, y, s)$ satisfy 
the first two equations in (\ref{FOC_barrQP}). Consequently, 
the duality gap is equal to the complementarity gap 
\begin{eqnarray}
   c^T x + \frac{1}{2} x^T Q x - (b^T y - \frac{1}{2} x^T Q x) = x^T s = n \mu, 
   \label{OptGap} 
\end{eqnarray}
and by reducing the barrier parameter $\mu$, IPM achieves convergence 
to optimality. The standard interior point algorithm applies the Newton 
method to (\ref{FOC_barrQP}), that is, computes the Newton direction and 
makes a step in this direction followed by a reduction of the barrier 
parameter $\mu$. The reduction of the barrier term is enforced by the use 
of the parameter $\sigma \in (0,1)$ and setting $\mu^{new} = \sigma \mu$. 
Due to the feasibility of $(x, y, s)$ the residual in (\ref{FOC_barrQP}) 
takes the following form 
\begin{eqnarray}
   (b-Ax, \, c-A^Ty-s+Qx, \, \sigma \mu e - XSe) = (0, 0, \xi). 
   \label{resNM}
\end{eqnarray}
Hence the Newton direction $(\Delta x, \Delta y, \Delta s)$ is obtained 
by solving the following system of linear equations 
\begin{eqnarray}
  \left[
  \begin{array}{rcl}
     A  &   0 & 0 \\
    -Q  & A^T & I \\
     S  &   0 & X
  \end{array}
  \right] \cdot \left[
  \begin{array}{l}
     \Delta x \\
     \Delta y \\
     \Delta s
  \end{array}
  \right] = \left[
  \begin{array}{c}
     0  \\
     0  \\
     \xi 
  \end{array}
  \right],
  \label{exactNM}
\end{eqnarray}
where $I$ denotes the identity matrix of dimension $n$.

In this paper we will analyse the method which allows the system 
(\ref{exactNM}) to be solved {\it inexactly}. To be precise, we will 
assume that all the iterates remain primal and dual feasible and 
an {\it inexact} Newton direction $(\Delta x, \Delta y, \Delta s)$ 
satisfies the following system of linear equations 
\begin{eqnarray}
  \left[
  \begin{array}{rcl}
     A  &   0 & 0 \\
    -Q  & A^T & I \\
     S  &   0 & X
  \end{array}
  \right] \cdot \left[
  \begin{array}{l}
     \Delta x \\
     \Delta y \\
     \Delta s
  \end{array}
  \right] = \left[
  \begin{array}{c}
     0  \\
     0  \\
     \xi + r 
  \end{array}
  \right]
  \label{inexactNM}
\end{eqnarray}
which admits an error $r$ in the third equation. Let us observe 
that any step in such a primal-dual inexact Newton direction 
preserves primal and dual feasibility. 

By using the positive semidefiniteness of $Q$ and exploiting the first 
two equations in the appropriate Newton system, it is easy to demonstrate 
that for both the exact (\ref{exactNM}) and the inexact (\ref{inexactNM}) 
Newton direction $(\Delta x, \Delta y, \Delta s)$, the following property 
holds: 
\begin{eqnarray}
  {\Delta x}^T \Delta s = {\Delta x}^T Q \, \Delta x \geq 0.
  \label{DxtDs}
\end{eqnarray}

The third equation in the Newton system plays a crucial role 
in the convergence analysis of an interior point algorithm. 
For the inexact Newton direction (\ref{inexactNM}) this equation 
takes the following form
\begin{eqnarray}
   S \Delta x + X \Delta s = \xi + r = \sigma \mu e - X S e + r, 
   \nonumber 
\end{eqnarray}
and by using $e^T e = n$ and $x^T s = n \mu$ we get 
\begin{eqnarray}
   s^T \Delta x + x^T \Delta s = \sigma \mu e^T e - x^T s + e^T r = (\sigma - 1) x^T s + e^T r. 
   \nonumber 
\end{eqnarray}
Hence the complementarity gap at the new point
$(x(\alpha), y(\alpha), s(\alpha)) = (x, y, s) + \alpha (\Delta x, \Delta y, \Delta s)$
becomes
\begin{eqnarray}
   x(\alpha)^T s(\alpha) & = & (x + \alpha \Delta x)^T (s + \alpha \Delta s) \nonumber \\
                         & = & x^T s + \alpha (s^T \Delta x + x^T \Delta s) 
                               + \alpha^2  \Delta x^T \Delta s \nonumber \\
                         & = & (1 - \alpha (1 - \sigma)) x^T s 
                               + \alpha e^T r + \alpha^2  \Delta x^T \Delta s 
   \label{bar_xts}
\end{eqnarray}
and the corresponding average complementarity gap is 
\begin{eqnarray}
   \mu(\alpha) = x(\alpha)^T s(\alpha) / n
               = (1 - \alpha ( 1 - \sigma )) \mu + \alpha e^T r /n + \alpha^2 {\Delta x}^T \Delta s /n.
   \label{NewGapAlpha}
\end{eqnarray}
Under the condition that the error term $e^T r$ and the second order 
term $\Delta x^T \Delta s$ are kept small enough in comparison 
with $\alpha (1 - \sigma) x^T s$, the complementarity gap at the new 
point is reduced compared with that at the previous iteration, 
thus guaranteeing progress of the algorithm. 

The convergence analysis of an interior point algorithm relies 
on imposing uniform progress in reducing the error in the complementarity 
conditions which technically is translated into a requirement that 
the error is small and bounded with ${\cal O}(\mu)$. To achieve it we will 
restrict the iterates to remain in the neighbourhood of the central 
path $\{ (x(\mu), y(\mu), s(\mu)), \, \mu > 0 \}$ and we will control 
the barrier reduction parameter $\sigma$ and the stepsize $\alpha$ 
so that $x(\alpha)^T s(\alpha)$ in (\ref{bar_xts}) is noticeably 
smaller than $x^T s$. 

We will consider two different ways of controlling the proximity 
to the central path and in both cases we will prove the convergence 
of the inexact interior point method and derive the worst-case complexity 
result. 
In both cases we will consider the {\it feasible} interior point 
algorithm and, hence, we will assume that all primal-dual iterates 
belong to the primal-dual strictly feasible set
${\cal F}^0 = \{ (x,y,s) \, | \, A x = b, \, A^T y + s - Q x = c, \, (x,s) > 0 \}$.
All iterates are confined to a neighbourhood of the central path which 
translates to a requirement that the error in the perturbed complementarity 
condition (\ref{FOC_barrQP}) is small. Depending on the norm used to measure 
this error, we will consider two different neighbourhoods: 
\begin{itemize}
\item 
a small neighbourhood induced by the use of the Euclidean norm for some $\theta \in (0,1)$
\begin{eqnarray}
   N_2(\theta) = \{ (x,y,s) \in {\cal F}^0 \, | \, \| X S e - \mu e \| \leq \theta \mu \}, 
   \label{N2hood}
\end{eqnarray}
which yields a short-step algorithm, and 
\item 
a symmetric neighbourhood induced by the use of the infinity norm for some $\gamma \in (0,1)$
\begin{eqnarray}
   N_S(\gamma) = \{ (x,y,s) \in {\cal F}^0 \, | 
                 \, \gamma \mu \leq x_j s_j \leq \frac{1}{\gamma} \mu, \, \forall j \},
   \label{NShood}
\end{eqnarray}
which yields a long-step algorithm.
\end{itemize}
The former has a theoretical importance as it leads to the algorithm 
with the best complexity result known to date. However, it is not 
an implementable method because it leads to poor performance in practice. 
Indeed, its behaviour reproduces the worst-case analysis. The latter 
neighbourhood has a practical meaning and leads to an efficient 
algorithm in practice. 

Following the general theory of the inexact Newton Method \cite{iter:DES,iter:Kelley},  
we will assume that the residual $r$ in (\ref{inexactNM}) satisfies 
\begin{eqnarray}
   \| r \|_p \leq \delta \| \xi \|_p, 
   \label{res}
\end{eqnarray}
for some $\delta \in (0,1)$ and an appropriate $p$-norm. 
According to the type of the neighbourhood used, (\ref{N2hood}) 
or (\ref{NShood}), this inequality will use either $p = 2$ 
or $p = \infty$, respectively. Further in the paper we will 
omit a subscript $2$ for the Euclidean norm unless an expression 
involves different norms at the same time and such an omission 
could lead to a confusion.

In the next section we will prove that the {\it feasible} interior 
point algorithm using an {\it inexact} \, Newton direction (\ref{inexactNM}) 
and applied to a convex quadratic program converges to an $\varepsilon$-accurate 
solution in ${\cal O}(\sqrt{n} \ln(1/\varepsilon))$ 
or ${\cal O}(n \ln(1/\varepsilon))$ iterations if it operates 
in the $N_2(\theta)$ or $ N_S(\gamma)$ neighbourhood, respectively. 
Our analysis will follow the general scheme used by Wright \cite{iter:Wright}.

\section{Worst-case complexity results}
\label{Complexity}

The analysis for two different neighbourhoods will share certain common 
features. The algorithm makes a step in the Newton direction obtained 
by solving (\ref{inexactNM}). When a step in the Newton direction 
$(\Delta x, \Delta y, \Delta s)$ is made, the new complementarity product 
for component $j$ is given by 
\begin{eqnarray}
   x_j(\alpha) \, s_j(\alpha) & = & (x_j + \alpha \Delta x_j) (s_j + \alpha \Delta s_j) \nonumber \\
                              & = & x_j s_j + \alpha (s_j \Delta x_j + x_j \Delta s_j) 
                                    + \alpha^2  \Delta x_j \Delta s_j. 
   \label{XjSj}
\end{eqnarray}
The third equation in (\ref{inexactNM}) is a local linearization 
of the complementarity condition and controls the middle term 
$s_j \Delta x_j + x_j \Delta s_j = \xi_j + r_j$ in the above equation. 
The error in the approximation of complementarity products is determined 
by the second-order term $\Delta x_j \Delta s_j$ in (\ref{XjSj}). 
We will provide a bound on these products, namely, we will bound 
the vector of the second-order error terms $\| \Delta X \Delta S e \|$. 

Having multiplied the third equation in the Newton system (\ref{inexactNM}) 
by $( X S )^{-1/2}$, we obtain
\begin{eqnarray}
   X^{-1/2} S^{1/2} \Delta x + X^{1/2} S^{-1/2} \Delta s = ( X S )^{-1/2} ( \xi + r ). 
   \label{3rdEq}
\end{eqnarray}
Defining $u = X^{-1/2} S^{1/2} \Delta x$ and $v = X^{1/2} S^{-1/2} \Delta s$ 
and using (\ref{DxtDs}) we obtain $u^T v = \Delta x^T \Delta s \geq 0$.
Let us partition all products $u_j v_j$ into positive and negative ones:
${\cal P} = \{ j \, | \, u_j v_j \geq 0 \}$ and ${\cal M} = \{ j \, | \, u_j v_j < 0 \}$
and observe that 
\begin{eqnarray}
    0 \leq u^T v 
      = \sum_{j \in {\cal P}}  u_j v_j  
      + \sum_{j \in {\cal M}}  u_j v_j 
      = \sum_{j \in {\cal P}} |u_j v_j| 
      - \sum_{j \in {\cal M}} |u_j v_j|. 
   \label{twoSums}
\end{eqnarray}
Next, let us write equation (\ref{3rdEq}) component-wise 
as $u_j + v_j = (x_j s_j)^{-1/2} (\xi_j + r_j)$ for every 
$j \in \{1, 2, \dots, n \}$ and take the sum of squared equations 
for components $j \in {\cal P}$: 
\begin{eqnarray*}
    0 \leq \sum_{j \in {\cal P}} (u_j + v_j)^2 
       =   \sum_{j \in {\cal P}} (u_j^2 + v_j^2) + 2 \sum_{j \in {\cal P}} u_j v_j  
       =   \sum_{j \in {\cal P}} (x_j s_j)^{-1} (\xi_j + r_j)^2, 
\end{eqnarray*}
to get 
\begin{eqnarray*}
     2 \sum_{j \in {\cal P}} | u_j v_j | 
     = 2 \sum_{j \in {\cal P}} u_j v_j 
     \leq \sum_{j \in {\cal P}} (x_j s_j)^{-1} (\xi_j + r_j)^2.
\end{eqnarray*}
Inequality (\ref{twoSums}) 
implies $\sum_{j \in {\cal M}} |u_j v_j| \leq \sum_{j \in {\cal P}} |u_j v_j|$
and, hence, we can write
\begin{eqnarray}
   \| \Delta X \Delta S e \|_1 
   & =    & \sum_{j \in {\cal P}} |u_j v_j| + \sum_{j \in {\cal M}} |u_j v_j| 
     \leq   2 \sum_{j \in {\cal P}} | u_j v_j | \nonumber \\
   & \leq & \sum_{j \in {\cal P}} (x_j s_j)^{-1} (\xi_j + r_j)^2 
     \leq   \sum_{j=1}^{n} (x_j s_j)^{-1} (\xi_j + r_j)^2.
   \label{DXDSeN1}
\end{eqnarray}

We will now consider two different algorithms: the short-step method 
in which the iterates are confined to $N_2(\theta)$ neighbourhood (\ref{N2hood}) 
and the long-step method in which the iterates are confined to $N_S(\gamma)$ 
neighbourhood (\ref{NShood}). The names ``short-step'' and ``long-step'' 
describe the steps to optimality and are related to the choice of barrier 
reduction parameter $\sigma$ (and should not be confused with the stepsizes 
taken in the Newton direction). 
In the short-step method, $\sigma$ is very close to one and therefore 
the algorithm makes only a short step to optimality while in the long step 
method, $\sigma$ is usually a small number satisfying $\sigma \ll 1$.

\subsection{Analysis of the short-step method}
\label{ShortSAnalysis}
In this section we will assume that $(x,y,s) \in N_{2}(\theta)$ 
for some $\theta \in (0,1)$ and inequality (\ref{res}) holds 
for the 2-norm: $\| r \|_2 \leq \delta \| \xi \|_2$.  
It is easy to deduce that since 
$\| X S e - \mu e \|_{\infty} \leq \| X S e - \mu e \|_2 \leq \theta \mu$, 
the complementarity products satisfy the following inequality 
\begin{eqnarray}
   (1 - \theta) \mu \leq x_j s_j \leq (1 + \theta) \mu \quad \forall j. 
   \label{XjSj_bound}
\end{eqnarray}
The barrier reduction parameter for the short-step algorithm is 
defined as $\sigma = 1 - \beta / \sqrt{n}$ for some $\beta \in (0,1)$. 
Such a definition implies that $(1 - \sigma)^2 n = \beta^2$ and therefore, 
using $e^T (X S e - \mu e) = 0$, the norm of the term $\xi = \sigma \mu e - XSe$ 
in (\ref{resNM}) satisfies
\begin{eqnarray}
   \| X S e - \sigma \mu e \|^2 
   & = & \| (X S e - \mu e ) + (1 - \sigma) \mu e \|^2 \nonumber \\ 
   & = & \| X S e - \mu e \|^2 
       + 2 (1 \! - \! \sigma) \mu e^T \! (X S e \! - \! \mu e) 
       + (1 \! - \! \sigma)^2 \mu^2 e^T e \nonumber \\ 
   & \leq & \theta^2 \mu^2 + (1 - \sigma)^2 n \mu^2 \nonumber \\
   & =    & (\theta^2 + \beta^2) \mu^2.
   \label{XIbound}
\end{eqnarray}
We are now ready to derive a bound on $\| \Delta X \Delta S e \|$.
\begin{lem}
\label{Distance}
Let $\theta \in (0,1)$. If $(x,y,s) \in N_{2}(\theta)$ then 
the inexact Newton direction $ (\Delta x, \Delta y, \Delta s)$ 
obtained by solving (\ref{inexactNM}) satisfies 
\begin{eqnarray}
   \| \Delta X \Delta S e \| \leq 
   \frac{(1 + \delta)^2 (\theta^2 +\beta^2)}{(1 - \theta)} \mu.
    \label{DXDSeN2}
\end{eqnarray}
\end{lem}

{\bf Proof:} 
Inequality (\ref{XjSj_bound}) provides a bound 
$\min_j \{x_j s_j\} \geq (1-\theta) \mu$. 
We use it to rewrite (\ref{DXDSeN1}): 
\begin{eqnarray}
   \| \Delta X \Delta S e \|_2 
    \leq \| \Delta X \Delta S e \|_1 
    \leq \frac{1}{\min_j \{x_j s_j\} } \sum_{j=1}^{n} (\xi_j + r_j)^2
    \leq \frac{1}{(1-\theta) \mu} \| \xi + r \|^2.
    \label{bnd2nrm}
\end{eqnarray}
Using (\ref{res}) (with $p = 2$) and (\ref{XIbound}) we write
\begin{eqnarray*}
    \| \xi + r \|^2 \leq ( \| \xi \| + \| r \| )^2 
                    \leq (1 + \delta)^2 \| \xi \|^2 
                    \leq (1 + \delta)^2 (\theta^2 + \beta^2) \mu^2, 
\end{eqnarray*}
and after substituting this expression into (\ref{bnd2nrm}) obtain 
the required inequality (\ref{DXDSeN2}).
\qed

Next we will show that for appropriately chosen constants $\theta, \beta$ 
and $\delta$, a full Newton step is feasible and the new iterate 
$(\bar x, \bar y, \bar s) = (x, y, s) + (\Delta x, \Delta y, \Delta s)$
also belongs to the $N_{2}(\theta)$ neighbourhood of the central path. 
We will prove an even stronger result which is that for any 
step $\alpha \in (0,1]$ in the Newton direction the following point 
\begin{eqnarray}
  (x(\alpha), y(\alpha), s(\alpha)) 
      = (x, y, s) + \alpha (\Delta x, \Delta y, \Delta s)
   \label{NewPointAlpha}
\end{eqnarray}
is primal-dual feasible and belongs to the $N_{2}(\theta)$ neighbourhood. 

Using (\ref{XjSj}) and (\ref{NewGapAlpha}) 
and $s_j \Delta x_j + x_j \Delta s_j = \sigma \mu - x_j s_j + r_j$, 
the deviation of the $j$th complementarity product from the average becomes 
\begin{eqnarray*}
   x_j(\alpha) s_j(\alpha) \! - \! \mu(\alpha) 
   & \!\! = \!\! & (x_j + \alpha \Delta x_j) (s_j + \alpha \Delta s_j) - \mu(\alpha) \\
   & \!\! = \!\! & x_j s_j + \alpha (s_j \Delta x_j + x_j \Delta s_j) 
                   + \alpha^2 \Delta x_j \Delta s_j - \mu(\alpha) \\
   & \!\! = \!\! & (1 \!\! - \! \alpha) x_j s_j \! + \! \alpha \sigma \mu 
                    \! + \! \alpha r_j \! + \! \alpha^2 \! \Delta x_j \Delta s_j 
       \! - \!\! (1 \!\! - \! \alpha) \mu \! - \! \alpha \sigma \mu 
                     \! - \! \alpha e^T r / n \! - \! \alpha^2 \! {\Delta x}^T \! \Delta s / n \\
   & \!\! = \!\! & (1 - \alpha) (x_j s_j - \mu) + \alpha (r_j - e^T r / n)
                   + \alpha^2 (\Delta x_j \Delta s_j - {\Delta x}^T \Delta s / n).
\end{eqnarray*}
Consequently, the proximity measure for the point $(x(\alpha), y(\alpha), s(\alpha))$ 
becomes
\begin{eqnarray}
   \| X(\alpha) S(\alpha)e \! - \! \mu(\alpha)e \| 
   \leq (1 \! - \! \alpha) \| X S e \! - \! \mu e \| 
        \! + \! \alpha \| r \! - \! \frac{ e^T r}{n} e \| 
        \! + \! \alpha^2 \| \Delta X \! \Delta S e \! - \! \frac{{\Delta x}^T \!\! \Delta s}{n} e \|.
   \label{NewProx}
\end{eqnarray}
The Lemma below sets three parameters: 
the proximity constant $\theta \in (0,1)$ in (\ref{N2hood}), 
$\beta \in (0,1)$ in the barrier reduction parameter 
$\sigma = 1 - \beta / \sqrt{n}$,  
and the level of error $\delta \in (0,1)$ allowed in the inexact 
Newton method (\ref{res}).

\begin{lem}
\label{NStepAlphaShortS}
Let $(x, y, s)$ be the current iterate in $N_{2}(\theta)$ neighbourhood 
and $(\Delta x, \Delta y, \Delta s)$ be the inexact Newton direction which 
solves equation system (\ref{inexactNM}). Let $\theta = \beta = 0.1$. 
If $\delta = 0.3$ then for any $\alpha \in (0,1]$ the new iterate 
(\ref{NewPointAlpha}) after a step $\alpha$ in this direction satisfies
\begin{eqnarray}
   \| X(\alpha) S(\alpha)e - \mu(\alpha)e \| 
   & \leq & \theta \mu(\alpha). 
   \label{NewProxAlpha}
\end{eqnarray}
\end{lem}

{\bf Proof:} 
Expanding the square and using $e^T e = n$ we write 
\begin{eqnarray}
   \| r \! - \! \frac{ e^T r}{n} e \|^2 \!\!
   & \!\! = \!\! & \!\! \| r \|^2 \! + \! \frac{1}{n^2} (e^T r)^2 e^T e 
                        \! - \! \frac{2}{n} (e^T r)^2 \nonumber \\
   & \!\! = \!\! & \!\! \| r \|^2 \! - \! \frac{1}{n} (e^T r)^2 
     \leq   \| r \|^2. 
   \label{RminusAvrR}
\end{eqnarray}
Similarly, expanding the square, using $e^T e = n$ again and 
$(\Delta X \Delta S e)^T e = {\Delta x}^T \Delta s$, we write
\begin{eqnarray}
   \| \Delta X \! \Delta S e \! - \! \frac{{\Delta x}^T \! \! \Delta s}{n} e \|^2  \!\!
   & \!\! = \!\! & \!\! \| \Delta X \Delta S e \|^2 + \frac{1}{n^2} ({\Delta x}^T \! \Delta s)^2 e^T \! e 
            - \frac{2 {\Delta x}^T \! \Delta s}{n} (\Delta X \Delta S e)^T e \nonumber \\
   & \!\! = \!\! & \!\! \| \Delta X \! \Delta S e \|^2 \! - \! \frac{1}{n} ({\Delta x}^T \!\! \Delta s)^2 
     \leq   \| \Delta X \! \Delta S e \|^2. 
   \label{dXdSminusAvr}
\end{eqnarray}
Next, using (\ref{NewProx}), the definition of $N_{2}(\theta)$ 
and inequalities (\ref{RminusAvrR}), (\ref{res}) and (\ref{dXdSminusAvr}), 
we write
\begin{eqnarray*}
   \| X(\alpha) S(\alpha)e \! - \! \mu(\alpha)e \| 
   & \leq & (1 \! - \! \alpha) \| X S e \! - \! \mu e \| 
            \! + \! \alpha \| r \! - \! \frac{ e^T r}{n} e \| 
            \! + \! \alpha^2 \| \Delta X \! \Delta S e 
            \! - \! \frac{{\Delta x}^T \!\! \Delta s}{n} e \| \nonumber \\
   & \leq & (1 \! - \! \alpha) \theta \mu 
            \! + \! \alpha \delta \| \xi \| 
            \! + \! \alpha^2 \| \Delta X \Delta S e \|. 
\end{eqnarray*}
Inequalities (\ref{XIbound}) and (\ref{DXDSeN2}) (Lemma~\ref{Distance}) 
provide bounds for the last two terms in the above inequality.
We use them to write 
\begin{eqnarray*}
   \| X(\alpha) S(\alpha)e \! - \! \mu(\alpha)e \| 
   & \leq & (1 \! - \! \alpha) \theta \mu 
            \! + \! \alpha \delta \sqrt{\theta^2 + \beta^2} \mu 
            \! + \! \alpha^2 \frac{(1 + \delta)^2 (\theta^2 +\beta^2)}{(1 - \theta)} \mu.
\end{eqnarray*}
The choice of $\theta = \beta = 0.1$ guarantees that $\sqrt{\theta^2 + \beta^2} \leq \sqrt{2} \theta$ 
and $\frac{(1 + \delta)^2 (\theta^2 +\beta^2)}{(1 - \theta)} = \frac{2(1 + \delta)^2}{9} \theta$
hence 
\begin{eqnarray*}
   \| X(\alpha) S(\alpha)e \! - \! \mu(\alpha)e \| 
   & \leq & (1 \! - \! \alpha) \theta \mu 
            \! + \! \sqrt{2} \alpha \delta \theta \mu 
            \! + \! \alpha^2 \frac{2 (1 + \delta)^2}{9} \theta \mu.
\end{eqnarray*}
Using equality (\ref{NewGapAlpha}) we observe that this lemma 
will be proved (inequality (\ref{NewProxAlpha}) will be satisfied) 
if the following holds
\begin{eqnarray*}
   (1 \! - \! \alpha) \theta \mu 
   \! + \! \sqrt{2} \alpha \delta \theta \mu 
   \! + \! \alpha^2 \frac{2 (1 + \delta)^2}{9} \theta \mu
   & \leq &
   \theta  \left( (1 - \alpha ( 1 - \sigma )) \mu + \alpha \frac{e^T r}{n} 
                + \alpha^2 \frac{ {\Delta x}^T \Delta s }{n} \right).
\end{eqnarray*}
We further simplify this inequality by removing the same terms 
present on both sides of it and then dividing it by $\alpha \theta$. 
Inequality (\ref{DxtDs}) guarantees that ${\Delta x}^T \Delta s$ is 
nonnegative, hence we conclude that (\ref{NewProxAlpha}) will be satisfied if
\begin{eqnarray*}
   \sqrt{2} \delta \mu \! + \! \alpha \frac{2 (1 + \delta)^2}{9} \mu
   & \leq &
   \sigma \mu + \frac{e^T r}{n}. 
\end{eqnarray*}
We observe that $|e^T r| \leq \| e \| \| r \| = \sqrt{n} \delta  \| \xi \|$ 
and hence, using (\ref{XIbound}), we write 
\begin{eqnarray}
   | \frac{e^T r}{n} |
   & \leq & \frac{\delta}{\sqrt{n}} \sqrt{\theta^2 + \beta^2} \mu 
   \ \leq \ \frac{\sqrt{2} \delta \theta}{\sqrt{n}} \mu.  
   \label{eTr}
\end{eqnarray}
Therefore to guarantee that (\ref{NewProxAlpha}) holds for any $\alpha \in (0,1]$,  
it suffices to choose $\delta$ such that 
\begin{eqnarray*}
   \sqrt{2} \delta \mu + \frac{2 (1 + \delta)^2}{9} \mu 
   & \leq &
   \sigma \mu - \frac{\sqrt{2} \delta \theta}{\sqrt{n}} \mu 
\end{eqnarray*}
and this simplifies to 
\begin{eqnarray*}
   \sqrt{2} \delta ( 1 \! + \! \frac{\theta}{\sqrt{n}} ) + \frac{2 (1 + \delta)^2}{9}
   & \leq &
   \sigma = 1 \! - \! \frac{\beta}{\sqrt{n}}. 
\end{eqnarray*}
The left hand side of this inequality is an increasing function 
of $\delta$ and we can easily check that the choice $\delta = 0.3$ gives 
$0.3 \sqrt{2} (1 \! + \! \theta / \sqrt{n}) + 3.38 / 9 \leq 1 \! - \! \beta / \sqrt{n}$ 
which holds for $\theta = \beta = 0.1$ and any $n \geq 2$. 
\qed

Lemma~\ref{NStepAlphaShortS} guarantees that for any $\alpha \in (0,1]$ 
the new iterate (\ref{NewPointAlpha}) also belongs to the $N_{2}(\theta)$ 
neighbourhood of the central path. We will set $\alpha = 1$ and 
take the full step in the Newton direction. 
Observe that the inexact Newton direction (\ref{inexactNM}) allows 
the error $r$ to appear only in the third equation which means that 
the direction $(\Delta x, \Delta y, \Delta s)$ preserves the feasibility 
of primal and dual equality constraints. The new iterate is defined as
$(\bar x, \bar y, \bar s) = (x, y, s) + (\Delta x, \Delta y, \Delta s)$ 
and setting $\alpha = 1$ in (\ref{NewGapAlpha}) gives 
\begin{eqnarray}
   \bar \mu = \mu(\alpha) = \sigma \mu + \frac{e^T r}{n} + \frac{{\Delta x}^T \Delta s}{n}.
   \label{barMu}
\end{eqnarray}

With $\theta = \beta = 0.1$ and $\delta = 0.3$ the right hand side 
term in inequality (\ref{DXDSeN2}) may be simplified to give 
$\| \Delta X \Delta S e \| \leq 0.38 \beta \mu$
and, using the Cauchy-Schwartz inequality, we get the bound 
\begin{eqnarray*}
   \frac{{\Delta x}^T \Delta s}{n}
   = \frac{(\Delta X \Delta S e)^T e}{n} 
   \leq \frac{1}{n} \| \Delta X \Delta S e \| \| e \| 
   \leq \frac{0.38}{\sqrt{n}} \beta \mu.  
\end{eqnarray*}
Using this inequality, our choice $\theta = \beta$, (\ref{eTr}) 
and $\sigma = 1 - \beta / \sqrt{n}$, 
we obtain the following bound on $\bar \mu$ in (\ref{barMu}) 
\begin{eqnarray}
   \bar \mu 
   \leq (1 - \frac{\beta}{\sqrt{n}}) \mu + \frac{2 \delta \beta}{\sqrt{n}} \mu 
                                         + \frac{0.38 \beta}{\sqrt{n}} \mu
    =   (1 - \frac{\eta}{\sqrt{n}}) \mu, 
   \label{NewGap}
\end{eqnarray}
where $\eta = \beta (1 - 2 \delta - 0.38) = 0.002$. 

We are now ready to state the complexity result for the inexact 
short-step feasible interior point method operating 
in a $N_{2}(0.1)$ neighbourhood.

\begin{teo}
\label{complThmShortS}
Given $\epsilon > 0$, suppose that a feasible starting point
$(x^0,y^0,s^0) \in N_{2}(0.1)$ satisfies
$(x^0)^T s^0 = n \mu^0, {\mbox {\rm { where }}} \mu^0 \leq 1/\epsilon^{\kappa}$,
for some positive constant $\kappa$. 
Then there exists an index $L$ with 
$L = {\cal O}(\sqrt{n} \, \ln (1/\epsilon) )$ such that 
$\mu^l \leq \epsilon, \ \forall l \geq L$.
\end{teo}

{\bf Proof:} 
is a straightforward application of Theorem 3.2 
in Wright \cite[Ch. 3]{iter:Wright}. 
\qed

\subsection{Analysis of the long-step method}
\label{LongSAnalysis}
In this section we will assume that $(x,y,s) \in N_S(\gamma)$
for some $\gamma \in (0,1)$ and inequality (\ref{res}) holds
for the infinity norm: $\| r \|_{\infty} \leq \delta \| \xi \|_{\infty}$.
We will ask for an aggressive reduction of the duality gap from one 
iteration to another and set the barrier reduction parameter 
$\sigma \in (0,1)$ to be significantly smaller than 1. 
Therefore it will not be possible in general to make the full step 
in the Newton direction. 

Using definition (\ref{NShood}) of the symmetric neighbourhood 
$N_S(\gamma)$ and observing that $1 / \gamma - 1 > 1 - \gamma$ 
we derive the following bound for the term $\xi = \sigma \mu e - XSe$
in (\ref{resNM})
\begin{eqnarray}
   \| X S e - \sigma \mu e \|_{\infty}
   &    = & \| (X S e - \mu e ) + (1 - \sigma) \mu e \|_{\infty} \nonumber \\
   & \leq & \| X S e - \mu e \|_{\infty} + (1 - \sigma) \mu \nonumber \\
   & \leq & \max \{ 1 - \gamma, \frac{1}{\gamma} - 1 \} \mu + (1 - \sigma) \mu \nonumber \\
   &    = & ( \frac{1}{\gamma} - \sigma ) \mu.
   \label{XIboundInf}
\end{eqnarray}
We are now ready to derive a bound on $\| \Delta X \Delta S e \|_{\infty}$.
\begin{lem}
\label{DistanceInf}
Let $\gamma \in (0,1)$. If $(x,y,s) \in N_{S}(\gamma)$ then
the inexact Newton direction $ (\Delta x, \Delta y, \Delta s)$
obtained by solving (\ref{inexactNM}) satisfies
\begin{eqnarray}
   \| \Delta X \Delta S e \|_{\infty} 
   \leq \| \Delta X \Delta S e \|_{1} 
   \leq n \frac{(1 + \delta)^2}{\gamma} ( \frac{1}{\gamma} - \sigma )^2 \mu
    \label{DXDSeNinfty}
\end{eqnarray}
and 
\begin{eqnarray}
   \Delta x_j \Delta s_j 
   \leq \frac{(1 + \delta)^2}{\gamma} ( \frac{1}{\gamma} - \sigma )^2 \mu, 
        \quad \forall j. 
    \label{xjsjNinfty}
\end{eqnarray}
\end{lem}

{\bf Proof:}
The definition of the $N_{S}(\gamma)$ neighbourhood provides a bound
$\min_j \{x_j s_j\} \geq \gamma \mu$.
We use it to rewrite (\ref{DXDSeN1}):
\begin{eqnarray}
   \| \Delta X \Delta S e \|_{\infty}
    \leq \| \Delta X \Delta S e \|_1
    \leq \frac{1}{\min_j \{x_j s_j\} } \sum_{j=1}^{n} (\xi_j + r_j)^2
    \leq \frac{1}{\gamma \mu} \| \xi + r \|_2^2.
    \label{bndInfnrm}
\end{eqnarray}
Using (\ref{res}) (for the infinity norm) and (\ref{XIboundInf}) we write
\begin{eqnarray*}
    \| \xi + r \|_2^2 \leq n \| \xi + r \|_{\infty}^2
                      \leq n (1 + \delta)^2 \| \xi \|_{\infty}^2
                      \leq n (1 + \delta)^2 ( \frac{1}{\gamma} - \sigma )^2 \mu^2, 
\end{eqnarray*}
and, after substituting this expression into (\ref{bndInfnrm}), 
obtain the required inequality (\ref{DXDSeNinfty}).
We observe that (\ref{DXDSeNinfty}) implies that 
\begin{eqnarray}
   - n \frac{(1 + \delta)^2}{\gamma} ( \frac{1}{\gamma} - \sigma )^2 \mu
       \leq \Delta x_j \Delta s_j 
       \leq n \frac{(1 + \delta)^2}{\gamma} ( \frac{1}{\gamma} - \sigma )^2 \mu, 
       \quad \forall j, 
    \nonumber
\end{eqnarray}
but we can obtain a tighter upper bound for this component-wise 
error term. For this we write equation (\ref{3rdEq}) for component $j$ 
and square both sides of it to get
\begin{eqnarray*}
    2 \Delta x_j \Delta s_j 
    \leq \frac{(\xi_j + r_j)^2}{x_j s_j}
    \leq \frac{(1 + \delta)^2 \| \xi \|_{\infty}^2}{\gamma \mu} 
    \leq \frac{(1 + \delta)^2}{\gamma} ( \frac{1}{\gamma} - \sigma )^2 \mu, 
    \quad \forall j, 
\end{eqnarray*}
which implies (\ref{xjsjNinfty}) and completes the proof.
\qed

Next we will show that for appropriately chosen constants $\sigma, \gamma$ 
and $\delta$, a (small) step $\alpha = {\cal O}(1/n)$ in the inexact 
Newton direction is feasible and the new iterate 
$(x(\alpha), y(\alpha), s(\alpha)) = (x, y, s) + \alpha (\Delta x, \Delta y, \Delta s)$ 
remains in the $N_{S}(\gamma)$ neighbourhood of the central path. 

The Lemma below provides conditions which have to be met by three 
parameters: the proximity constant $\gamma \in (0,1)$ in (\ref{NShood}),
the barrier reduction parameter $\sigma \in (0,1)$, 
and the level of error $\delta \in (0,1)$ allowed in the inexact
Newton method (\ref{res}).

\begin{lem}
\label{NStepAlphaLongS}
Let $(x, y, s)$ be the current iterate in the $N_{S}(\gamma)$ neighbourhood
and $(\Delta x, \Delta y, \Delta s)$ be the inexact Newton direction 
which solves equation system (\ref{inexactNM}). 
If the stepsize $\alpha \in (0,1]$ satisfies the following conditions: 
\begin{eqnarray}
   \alpha (\gamma + n) \frac{(1 + \delta)^2}{\gamma} ( \frac{1}{\gamma} - \sigma )^2     
      & \leq & \sigma (1 - \gamma) 
      - \delta ( 1 + \gamma ) (\frac{1}{\gamma} - \sigma) \label{cnd1} \\
   \delta ( 1 + \frac{1}{\gamma} ) (\frac{1}{\gamma} - \sigma) 
   + \alpha \frac{(1 + \delta)^2}{\gamma} ( \frac{1}{\gamma} - \sigma )^2 
      & \leq & (\frac{1}{\gamma} - 1) \sigma  \label{cnd2}
\end{eqnarray}
then the new iterate $(x(\alpha), y(\alpha), s(\alpha))$  
belongs to the $N_{S}(\gamma)$ neighbourhood, that is: 
\begin{eqnarray}
   \gamma \mu(\alpha) 
       \leq x_j(\alpha) s_j(\alpha) 
       \leq \frac{1}{\gamma} \mu(\alpha), \quad \forall j. 
   \label{inNSgamma}
\end{eqnarray}
\end{lem}

{\bf Proof:}
Using the average complementarity gap (\ref{NewGapAlpha}) at the new 
iterate $(x(\alpha), y(\alpha), s(\alpha))$ together with expression 
(\ref{XjSj}) and the third equation in (\ref{inexactNM}), we deduce 
that the left inequality in (\ref{inNSgamma}) will hold for any $j$ if 
\begin{eqnarray}
   \gamma \left( (1 - \alpha) \mu + \alpha \sigma \mu + \alpha \frac{e^T r}{n} 
            + \alpha^2 \frac{{\Delta x}^T \Delta s}{n} \right)
   \leq 
   (1 - \alpha) x_j s_j + \alpha \sigma \mu + \alpha r_j 
                        + \alpha^2 \Delta x_j \Delta s_j. 
   \nonumber
\end{eqnarray}
Since $(x, y, s) \in N_{S}(\gamma)$ we know that 
$\gamma (1 - \alpha) \mu \leq (1 - \alpha) x_j s_j$ and therefore 
the above inequality will hold if we satisfy a tighter version of it: 
\begin{eqnarray}
   \gamma \left( (1 - \alpha) \mu + \alpha \sigma \mu + \alpha \frac{e^T r}{n}
            + \alpha^2 \frac{{\Delta x}^T \Delta s}{n} \right)
   \leq
   \gamma (1 - \alpha) \mu + \alpha \sigma \mu + \alpha r_j
                           + \alpha^2 \Delta x_j \Delta s_j.
   \nonumber
\end{eqnarray}
After removing identical terms from both sides and dividing 
both sides by $\alpha$, we conclude that the inequality will hold if 
\begin{eqnarray}
   \alpha \left( \gamma \frac{{\Delta x}^T \Delta s}{n} - \Delta x_j \Delta s_j \right) 
   \leq
   \sigma (1 - \gamma) \mu + r_j - \gamma \frac{e^T r}{n}. 
   \label{ineq101}
\end{eqnarray}
Using Lemma~\ref{DistanceInf} we deduce that 
\begin{eqnarray}
   \frac{{\Delta x}^T \Delta s}{n}
   = \frac{(\Delta X \Delta S e)^T e}{n}
   \leq \frac{1}{n} \| \Delta X \Delta S e \|_1 \| e \|_{\infty}
   \leq \frac{(1 + \delta)^2}{\gamma} ( \frac{1}{\gamma} - \sigma )^2 \mu 
   \label{ineq102}
\end{eqnarray}
hence 
\begin{eqnarray}
   \gamma \frac{{\Delta x}^T \Delta s}{n} - \Delta x_j \Delta s_j 
   \leq (\gamma + n) \frac{(1 + \delta)^2}{\gamma} ( \frac{1}{\gamma} - \sigma )^2 \mu. 
   \label{ineq103}
\end{eqnarray}
Using (\ref{res}) (for the infinity norm) and (\ref{XIboundInf}) we deduce that 
$\| r \|_{\infty} \leq \delta \| \xi \|_{\infty} \leq \delta (\frac{1}{\gamma} - \sigma) \mu$ 
and 
\begin{eqnarray}
   | e^T r | \leq \| r \|_1 
   \leq n \| r \|_{\infty} \leq n \delta (\frac{1}{\gamma} - \sigma) \mu
   \label{ineq104}
\end{eqnarray}
hence 
\begin{eqnarray}
   | r_j - \gamma \frac{e^T r}{n} | 
   \leq | r_j | + \gamma | \frac{e^T r}{n} | 
   \leq \delta ( 1 + \gamma ) (\frac{1}{\gamma} - \sigma) \mu. 
   \label{ineq105}
\end{eqnarray}
%
The inequalities (\ref{ineq103}) and (\ref{ineq105}) allow us 
to determine the most adverse conditions in (\ref{ineq101}), 
namely when its left-hand-side is the largest possible 
and the right-hand-side is the smallest possible: 
\begin{eqnarray}
   \alpha (\gamma + n) \frac{(1 + \delta)^2}{\gamma} ( \frac{1}{\gamma} - \sigma )^2 \mu 
      & \leq & \sigma (1 - \gamma) \mu 
      - \delta ( 1 + \gamma ) (\frac{1}{\gamma} - \sigma) \mu. 
   \nonumber 
\end{eqnarray}
We then conclude that (\ref{cnd1}) implies (\ref{ineq101}) 
and therefore the left inequality in (\ref{inNSgamma}) holds. 
This completes the first part of the proof. 

The second part deals with the right inequality in (\ref{inNSgamma}). 
Again, using (\ref{NewGapAlpha}), (\ref{XjSj}) and the third equation 
in (\ref{inexactNM}), we write the required inequality which should 
be satisfied for any $j \in \{ 1, 2, \dots, n \}$ 
\begin{eqnarray}
   (1 - \alpha) x_j s_j + \alpha \sigma \mu + \alpha r_j
                        + \alpha^2 \Delta x_j \Delta s_j
   \leq
   \frac{1}{\gamma} \left( (1 - \alpha) \mu + \alpha \sigma \mu 
       + \alpha \frac{e^T r}{n} + \alpha^2 \frac{{\Delta x}^T \Delta s}{n} \right), 
   \nonumber
\end{eqnarray}
and determine the condition under which it holds. 
We use similar arguments as in the earlier part of the proof: 
for example, $(x, y, s) \in N_{S}(\gamma)$ implies 
$(1 - \alpha) x_j s_j \leq \frac{1}{\gamma} (1 - \alpha) \mu$. 
Hence we simplify this inequality by removing identical terms 
and then divide both sides of it by $\alpha$. 
The inequality we need to satisfy becomes 
\begin{eqnarray}
   r_j - \frac{e^T r}{\gamma n}
   + \alpha \left( \Delta x_j \Delta s_j - \frac{{\Delta x}^T \Delta s}{\gamma n} \right)
   \leq
   (\frac{1}{\gamma} - 1) \sigma \mu.
   \label{ineq110}
\end{eqnarray}
Using (\ref{xjsjNinfty}) and (\ref{DxtDs}) we write 
\begin{eqnarray}
   \Delta x_j \Delta s_j - \frac{{\Delta x}^T \Delta s}{\gamma n} 
   \leq \frac{(1 + \delta)^2}{\gamma} ( \frac{1}{\gamma} - \sigma )^2 \mu 
   \label{ineq111}
\end{eqnarray}
and using (\ref{res}) and (\ref{XIboundInf})
and similar arguments to those which
led to (\ref{ineq105}) we write
\begin{eqnarray}
   | r_j - \frac{e^T r}{\gamma n} |
   \leq | r_j | + | \frac{e^T r}{\gamma n} |
   \leq \delta ( 1 + \frac{1}{\gamma} ) (\frac{1}{\gamma} - \sigma) \mu.
   \label{ineq112}
\end{eqnarray}
By (\ref{ineq111}) and (\ref{ineq112}), to satisfy (\ref{ineq110}) 
we need to choose $\alpha$ such that: 
\begin{eqnarray}
   \delta ( 1 + \frac{1}{\gamma} ) (\frac{1}{\gamma} - \sigma) \mu
   + \alpha \frac{(1 + \delta)^2}{\gamma} ( \frac{1}{\gamma} - \sigma )^2 \mu
   \leq (\frac{1}{\gamma} - 1) \sigma \mu 
   \nonumber 
\end{eqnarray}
which simplifies to (\ref{cnd2}) and completes the second part of the proof.
\qed

Lemma~\ref{NStepAlphaLongS} provides conditions which the stepsize 
$\alpha \in (0,1]$ needs to satisfy so that the new iterate 
remains in the $N_{S}(\gamma)$ neighbourhood of the central path. 
We still need to demonstrate that after a step is made 
a sufficient reduction of duality gap is achieved. 

\begin{lem}
\label{gapRedLongS}
Let $(x, y, s) \in N_{S}(\gamma)$ be given and let 
$(\Delta x, \Delta y, \Delta s)$ be the inexact Newton direction 
which solves equation system (\ref{inexactNM}). 
If the stepsize $\alpha \in (0,1]$ satisfies the inequality 
\begin{eqnarray}
   \sigma + \delta (\frac{1}{\gamma} - \sigma)  
          + \alpha \frac{(1 + \delta)^2}{\gamma} ( \frac{1}{\gamma} - \sigma )^2 
   \leq 0.9,  
   \label{cnd3} 
\end{eqnarray}
then the duality gap at the new iterate $(x(\alpha), y(\alpha), s(\alpha))$  
satisfies: 
\begin{eqnarray}
   \mu(\alpha) \leq (1 - 0.1 \alpha) \mu. 
   \label{gapRed}
\end{eqnarray}
\end{lem}

{\bf Proof:}
By substituting (\ref{NewGapAlpha}) into (\ref{gapRed}), 
cancelling similar terms and dividing the resulting inequality 
by $\alpha$, we replace (\ref{gapRed}) with a new condition 
that the stepsize $\alpha$ has to satisfy: 
\begin{eqnarray}
   \sigma \mu + \frac{e^T r}{n} + \alpha \frac{{\Delta x}^T \Delta s}{n} 
   \leq 0.9 \mu. 
   \nonumber
\end{eqnarray}
Using the bounds (\ref{ineq104}) and (\ref{ineq102}) derived earlier 
we conclude that the above inequality will hold if 
\begin{eqnarray}
   \sigma \mu + \delta (\frac{1}{\gamma} - \sigma) \mu 
              + \alpha \frac{(1 + \delta)^2}{\gamma} ( \frac{1}{\gamma} - \sigma )^2 \mu 
   \leq 0.9 \mu, 
   \nonumber
\end{eqnarray}
which is equivalent to (\ref{cnd3}). 
\qed

It remains to consider the three conditions (\ref{cnd1}), (\ref{cnd2}) 
and (\ref{cnd3}) and to demonstrate that an appropriate choice 
of parameters $\gamma, \sigma$ and $\delta$ guarantees that 
all these conditions hold for some $\alpha = {\cal O}(\frac{1}{n})$. 

We set the proximity constant $\gamma = 0.5$ in (\ref{NShood}), 
the barrier reduction parameter $\sigma = 0.5$ and 
$\delta = 0.05$ as the level of error allowed in the inexact
Newton method (\ref{res}). Indeed, with these parameter settings 
%
$
   \frac{(1 + \delta)^2}{\gamma} ( \frac{1}{\gamma} - \sigma )^2 
   = 4.96125, 
$
%
and we verify that all three conditions (\ref{cnd1}), (\ref{cnd2}) 
and (\ref{cnd3}) are satisfied by $\hat \alpha = \frac{1}{50n}$
for any $n \geq 2$. Substituting such an $\hat \alpha$ into 
(\ref{gapRed}) gives 
\begin{eqnarray}
   \bar \mu = \mu(\hat \alpha) \leq (1 - \frac{\eta}{n}) \mu, 
   \nonumber 
\end{eqnarray}
where $\eta = 0.002$, and allows us to conclude this section 
with the following complexity result for the long-step inexact 
feasible interior point method operating in a $N_{S}(0.5)$ 
neighbourhood.

\begin{teo}
\label{complThmLongS}
Given $\epsilon > 0$, suppose that a feasible starting point
$(x^0,y^0,s^0) \in N_{S}(0.5)$ satisfies
$(x^0)^T s^0 = n \mu^0, {\mbox {\rm { where }}} \mu^0 \leq 1/\epsilon^{\kappa}$,
for some positive constant $\kappa$.
Then there exists an index $L$ with
$L = {\cal O}(n \, \ln (1/\epsilon) )$ such that
$\mu^l \leq \epsilon, \ \forall l \geq L$.
\end{teo}

{\bf Proof:}
is a straightforward application of Theorem 3.2
in Wright \cite[Ch. 3]{iter:Wright}.
\qed

\section{Practical aspects of the inexact IPM}
\label{PracIPM} 

Linear system (\ref{inexactNM}) is a modification of (\ref{exactNM}) 
which admits the error $r$ in its third equation. The analysis 
in the previous section provides guarantees that when the Newton 
equation system (\ref{inexactNM}) is solved inexactly and error $r$ 
satisfies condition (\ref{res}) the interior point algorithms retain 
their good complexity results. In this section we will briefly 
comment on a practical way of computing the Newton direction 
which meets condition (\ref{res}). 

Let us observe that after eliminating $\Delta s$, (\ref{exactNM}) 
is reduced to the augmented form 
\begin{eqnarray}
   \left[
   \begin{array}{cc}
     - Q - \Theta^{-1} & A^T \\
       A               & 0
   \end{array}
   \right]
   \left[
   \begin{array}{c}
      \Delta x \\
      \Delta y
   \end{array}
   \right]
   =
   \left[
   \begin{array}{c}
      - X^{-1} \xi \\
      0
   \end{array}
   \right],
   \label{AugSys}
\end{eqnarray}
where $\Theta = X S^{-1} \in {\cal R}^{n \times n}$.
Any IPM has to solve at least one such system at each 
iteration \cite{iter:dAdSdS,iter:G-ipmXXV}. 
Numerous attempts have been made during the last decade
to employ an iterative method for this task. 
Iterative methods for linear algebra are particularly attractive 
if they can be used to find only an approximate solution 
of the linear system, that is, when their run can be truncated  
to merely a few iterations. It is common to interrupt 
the iterative process once the required (loose) accuracy 
of the solution is obtained. Clearly such a solution 
is inexact and in the context of (\ref{AugSys}) this 
translates to dealing with an inexact Newton direction 
$(\Delta \tilde x, \Delta \tilde y)$ which satisfies 
\begin{eqnarray}
   \left[
   \begin{array}{cc}
     - Q - \Theta^{-1} & A^T \\
       A               & 0
   \end{array}
   \right]
   \left[
   \begin{array}{c}
      \Delta \tilde x \\
      \Delta \tilde y
   \end{array}
   \right]
   =
   \left[
   \begin{array}{c}
      - X^{-1} \xi + r_x \\
      r_y 
   \end{array}
   \right],
   \label{inexactAugSys}
\end{eqnarray}
where the errors $r_x \in {\cal R}^n$ and $r_y \in {\cal R}^m$ 
determine the level of inexactness. 

The analysis presented in this paper applies to the situation
when $r_y = 0$. The other error, $r_x$ may take a nonzero value and 
indeed, (\ref{inexactAugSys}) becomess equivalent to (\ref{inexactNM})
if $r_y = 0$ and
\begin{eqnarray}
   - X^{-1} \xi + r_x = - X^{-1} (\xi + r).
\end{eqnarray}
This equation combined with condition (\ref{res}) determines 
practical stopping criteria set for an iterative solution method 
applied to (\ref{AugSys}): 
\begin{eqnarray}
   r_y = 0 \quad {\mbox { and }} \quad 
   \| r \| = \| X r_x \| \leq \delta \| \xi \|.
\end{eqnarray}
Interestingly, there exist classes of iterative methods which can 
meet the above stopping criteria. They belong to a broad collection 
of iterative methods for saddle point problems \cite{iter:BenziGL} 
and rely on specially designed indefinite 
preconditioners \cite{iter:LV-Indef,iter:RS-Krylov}.

\section{Conclusions}
\label{Conclusions} 

The analysis presented in this paper provides the proofs 
of ${\cal O}(\sqrt{n} \log(1/\varepsilon))$ and 
${\cal O}(n \log(1/\varepsilon))$ iteration complexity 
of the, respectively, short-step and long-step inexact feasible 
primal-dual algorithms for quadratic programming. The analysis 
allows for considerable relative errors in the Newton direction. 
Indeed, $\delta$ in (\ref{res}) may take values $0.3$ 
and $0.05$ for the short-step and long-step algorithms, 
respectively. This shows, somewhat surprisingly, that 
the inexactness in the solution of the Newton equation 
system (\ref{inexactNM}) may be quite considerable without 
adversly affecting the best known worst-case iteration 
complexity results of these algorithms. It is an encouraging 
result for researchers who design preconditioners for 
iterative methods and wish to apply them to solve the reduced 
Newton equation systems arising in the context of interior 
point methods.

{
{\bf Acknowledgement} \\
The author is grateful to 
Kimonas Fountoulakis, Pablo Gonz\'alez-Brevis and Julian Hall 
for reading earlier versions of this paper and useful suggestions 
which led to an improvement of the presentation. 
}

%

\end{document}